\documentclass[12pt]{article}
\usepackage{amsfonts}
\usepackage{mathrsfs}

\usepackage{amsmath,amssymb,amsfonts,theorem,makeidx,latexsym,epsfig,color,subfigure,algorithm,algorithmic}\usepackage[numbers,sort&compress]{natbib}

\newtheorem{defn}{Definition}[section]

\newtheorem{lem}[defn]{Lemma}

{\theorembodyfont{\rmfamily}

}

\newtheorem{thm}[defn]{Theorem}

\newtheorem{obs}[defn]{Observation}

\newtheorem{cor}[defn]{Corollary}

\numberwithin{equation}{section}

\def\bp{{\noindent\bf Proof. \ }}

\def\ep{\hfill$\square$\par\bigskip}

\title{A note on semitotal domination in graphs \footnote{The research is supported by NSFC (No. 11301440),
Natural Science Foundation of Fujian Province (CN)(2015J05017)}}

\author{{\large Wei Zhuang \thanks{Corresponding author; E-mail: zhuangweixmu@163.com}} \\
{\it \normalsize School of Applied Mathematics},
\\{\it \normalsize  Xiamen University of Technology, Xiamen 361024, P.R.China }}

\date{}

\begin{document}

\maketitle

\begin{abstract} A set $S$ of vertices in $G$ is a semitotal
dominating set of $G$ if it is a dominating set of $G$ and every
vertex in $S$ is within distance $2$ of another vertex of $S$. The
\emph{semitotal domination number}, $\gamma_{t2}(G)$, is the minimum
cardinality of a semitotal dominating set of $G$. The \emph{semitotal domination multisubdivision number} of a
graph $G$, $msd_{\gamma_{t2}}(G)$, is the minimum positive integer $k$ such that there exists an edge which must be subdivided $k$ times to increase the semitotal
domination number of $G$. In this paper, we show that $msd_{\gamma_{t2}}(G)\leq 3$ for any graph $G$ of order at least
$3$, we also determine the semitotal domination multisubdivision number for some
classes of graphs and characterize trees $T$ with $msd_{\gamma_{t2}}(T)=3$. On the other hand, we know that $\gamma_{t2}(G)$ is a parameter that is squeezed between domination number, $\gamma(G)$ and total domination number, $\gamma_t(G)$, so for any tree $T$, we investigate the ratios $\frac{\gamma_{t2}(T)}{\gamma(T)}$ and $\frac{\gamma_t(T)}{\gamma_{t2}(T)}$, and present the constructive characterizations of the families of trees achieving the
upper bounds.
\end{abstract}

\begin{minipage}{150mm}
{\bf Keywords}\ {domination number, semitotal domination number, semitotal domination multisubdivision number}\\

{\bf AMS Subject Classification:}\ {05C05, 05C69}\\

\end{minipage}

\section{Introduction}
Let $G=(V, E)$ be a simple graph without isolated vertices, and let
$v$ be a vertex in $G$. The \emph{open neighborhood} of $v$ is
$N(v)=\{u\in V|uv\in E\}$ and the \emph{degree} of $v$ is
$d(v)=|N(v)|$. For two vertices $u$ and $v$ in a connected graph
$G$, the \emph{distance} $d(u, v)$ between $u$ and $v$ is the length
of a shortest $(u, v)$-path in $G$. The maximum distance among all
pairs of vertices of $G$ is the \emph{diameter} of a graph $G$ which
is denoted by $diam(G)$. A \emph{leaf} of $G$ is a vertex of degree
$1$, and a \emph{support vertex} of $G$ is a vertex adjacent to a
leaf. A vertex $v$ is called \emph{universal} if $d(v)=n-1$.

A set $S$ of vertices of a graph $G$ is called a \emph{dominating
set} (respectively, \emph{total dominating set}) of $G$ if every
vertex in $V(G)\setminus S$ (respectively, $V(G)$) is adjacent to at
least one vertex in $S$. The \emph{domination number} (respectively,
\emph{total domination number}) of $G$, denoted by $\gamma(G)$
(respectively, $\gamma_t(G)$), is the minimum cardinality of a
dominating set (respectively, total dominating set) of $G$.

The concept of semitotal domination in graphs was introduced by
Goddard et al.\cite{Goddard}. A set $S$ of vertices in $G$ is a
\emph{semitotal dominating set} of $G$ if it is a dominating set of
$G$ and every vertex in $S$ is within distance $2$ of another vertex
of $S$. The \emph{semitotal domination number}, $\gamma_{t2}(G)$, is
the minimum cardinality of a semitotal dominating set of $G$. We
observe that $\gamma(G)\leq \gamma_{t2}(G)\leq \gamma_t(G)$. A
semitotal dominating set (respectively, dominating set, total
dominating set) of $G$ of cardinality $\gamma_{t2}(G)$
(respectively, $\gamma(G)$, $\gamma_t(G)$) is called a
\emph{$\gamma_{t2}(G)$-set} (respectively, \emph{$\gamma(G)$-set},
\emph{$\gamma_t(G)$-set}).

The \emph{semitotal domination multisubdivision number} of a
graph $G$, $msd_{\gamma_{t2}}(G)$, is the minimum positive integer $k$ such that there exists
an edge which must be subdivided $k$ times to increase the semitotal
domination number of $G$. Similar definitions exist for the
domination multisubdivision number and total
domination multisubdivision number, which are introduce in \cite{Dettlaff} and \cite{Alaminos} respectively.

One of the purposes of our paper is to
initialize the study of the seimitotal domination multisubdivision
number. We show that $msd_{\gamma_{t2}}(G)\leq 3$ for any connected graph $G$ of order at least
$3$, we also determine the semitotal domination multisubdivision number for some
classes of graphs and characterize trees $T$ with $msd_{\gamma_{t2}}(T)=3$.

On the other hand, an area of research in domination of graphs that has received
considerable attention is the study of the ratio between two domination
parameters (some related results can be referred to \cite{Cyman1, Cyman2,
Chellali, Krzywkowski, Henning3, Zhu, Dayila}). Combining this with the fact that $\gamma(G)\leq \gamma_{t2}(G)\leq \gamma_t(G)$, it is natural to consider the ratios $\frac{\gamma_{t2}(T)}{\gamma(T)}$ and $\frac{\gamma_t(T)}{\gamma_{t2}(T)}$. We also present the constructive characterizations of the families of trees achieving the
upper bounds.

\section{Results and bounds for the semitotal domination multisubdivision
number}

\subsection{Preliminary results}

In this section, we start with some basic results of semitotal domination multisubdivision number. From the definition of semitotal domination multisubdivision number, the following conclusions are trivial.

\begin{obs}
For a complete graph $K_n$ and a wheel $W_n$, $n\geq 3$, we have that $msd_{\gamma_{t2}}(K_n)=msd_{\gamma_{t2}}(W_n)=
3$.
\end{obs}

\begin{obs}
For a cycle $C_n$ and a path $P_n$, $n\geq 3$, we have that
$$msd_{\gamma_{t2}}(P_n)=msd_{\gamma_{t2}}(C_n)= \left\{
          \begin{array}{ll}
            1, & \hbox{if $n\equiv 0, 2 ($mod $5)$;} \\
            2, & \hbox{if $n\equiv 1, 4 ($mod $5)$;} \\
            3, & \hbox{if $n\equiv 3 ($mod $5)$}
          \end{array}
        \right.
        $$
\end{obs}

\begin{obs}
For a complete bipartite graph $K_{p, q}$, $p\leq q$, we have that
$$msd_{\gamma_{t2}}(K_{p, q})= \left\{
          \begin{array}{ll}
            4, & \hbox{if $p=q=1$;} \\
            3, & \hbox{if $p=1$ and $q>1$;} \\
            2, & \hbox{if $p\geq 2$}
          \end{array}
        \right.
        $$
\end{obs}

The main result of this section is the next theorem.

\begin{thm}
For a connected graph $G$ of order at least $3$, $msd_{\gamma_{t2}}(G)\leq 3$.
\end{thm}

\bp
We take an edge of $G$, say $uv$, and subdivide it with vertices $u_1, u_2, u_3$. Denote the resulting graph by $G'$, let $D'$ be a $\gamma_{t2}$-set of $G'$. Clearly, $|\{u_1, u_2, u_3\}\cap D'|\leq 2$. We consider the next two cases.

{\flushleft\textbf{Case 1.}}\quad $|\{u_1, u_2, u_3\}\cap D'|=1$.

If $\{u_1, u_2, u_3\}\cap D'=\{u_2\}$, then either $u$ or $v$ belongs to $D'$, say $u$. Moreover, $(N_{G'}[v]\setminus \{u_3\})\cap D' \neq \emptyset$. We set $D=D'\setminus \{u_2\}$.

If $\{u_1, u_2, u_3\}\cap D'=\{u_1\}$, then $(N_{G'}[u]\setminus \{u_1\})\cap D' \neq \emptyset$. Moreover, $v$ belongs to $D'$. We set $D=D'\setminus \{u_1\}$.

In either case, the set $D$ is a semitotal dominating set of $G$, and so $\gamma_{t2}(G)\leq \gamma_{t2}(G')-1$.

{\flushleft\textbf{Case 2.}}\quad $|\{u_1, u_2, u_3\}\cap D'|=2$.

Let $D''=(D'\setminus \{u_2\})\cup \{u_1, u_3\}$, and clearly, $D''$ is also a $\gamma_{t2}$-set of $G'$. Since $|G|\geq 3$, either $u$ or $v$ has degree at least two, say $u$. Assume that $w$ is a neighbor of $u$ in $G$ other than $v$. Note that $N_{G'}[w]\cap D'' \neq \emptyset$, set
$D'''=(D''\setminus \{u_1, u_3\})\cup \{u\}$ when $u\not \in D''$, and $D'''=(D''\setminus \{u_1, u_3\})\cup \{w\}$ when $u\in D''$.
Clearly, $D'''$ is a semitotal dominating set of $G$, and then $\gamma_{t2}(G)\leq \gamma_{t2}(G')-1$.
\ep

\begin{cor}
If there is a universal vertex in a graph $G$ of order at least $3$, then $msd_{\gamma_{t2}}(G)=3$.
\end{cor}

Moreover, if $T$ is a tree, then we have the following results.

\begin{obs}
$(1)$ Let $T$ be a tree, $u, v$ be two adjacent support vertices. If
either $u$ or $v$ has degree two, then $msd_{\gamma_{t2}}(T)\leq 2$.

$(2)$ Let $T$ be a tree, $u, v$ be two support vertices at distance
two apart. If either $u$ or $v$ has degree two, then
$msd_{\gamma_{t2}}(T)\leq 2$.
\end{obs}

\subsection{Trees with semitotal domination multisubdivision number equal to 3}

From Theorem~2.4, we know that for a tree $T$ of order at least $3$, $msd_{\gamma_{t2}}(T)\leq 3$.
Trees are classified as Class~1, Class~2 and Class~3 depending on
whether their semitotal domination multisubdivision number is 1, 2
or 3, respectively. In the following, we are ready to provide a
constructive characterization of trees in Class~3.

For our purposes, we
define a \emph{labeling} of a tree $T$ as a partition $S=(S_A, S_B,
S_C)$ of $V(T)$ (This idea of labeling the vertices is
introduced in \cite{Dorfling}). We will refer to the pair $(T, S)$
as a \emph{labeled tree}. The label or \emph{status} of a vertex
$v$, denoted sta$(v)$, is the letter $x\in \{A, B, C\}$ such
that $v\in S_x$.

$\\$

\begin{center}
  \includegraphics[width=2.5in]{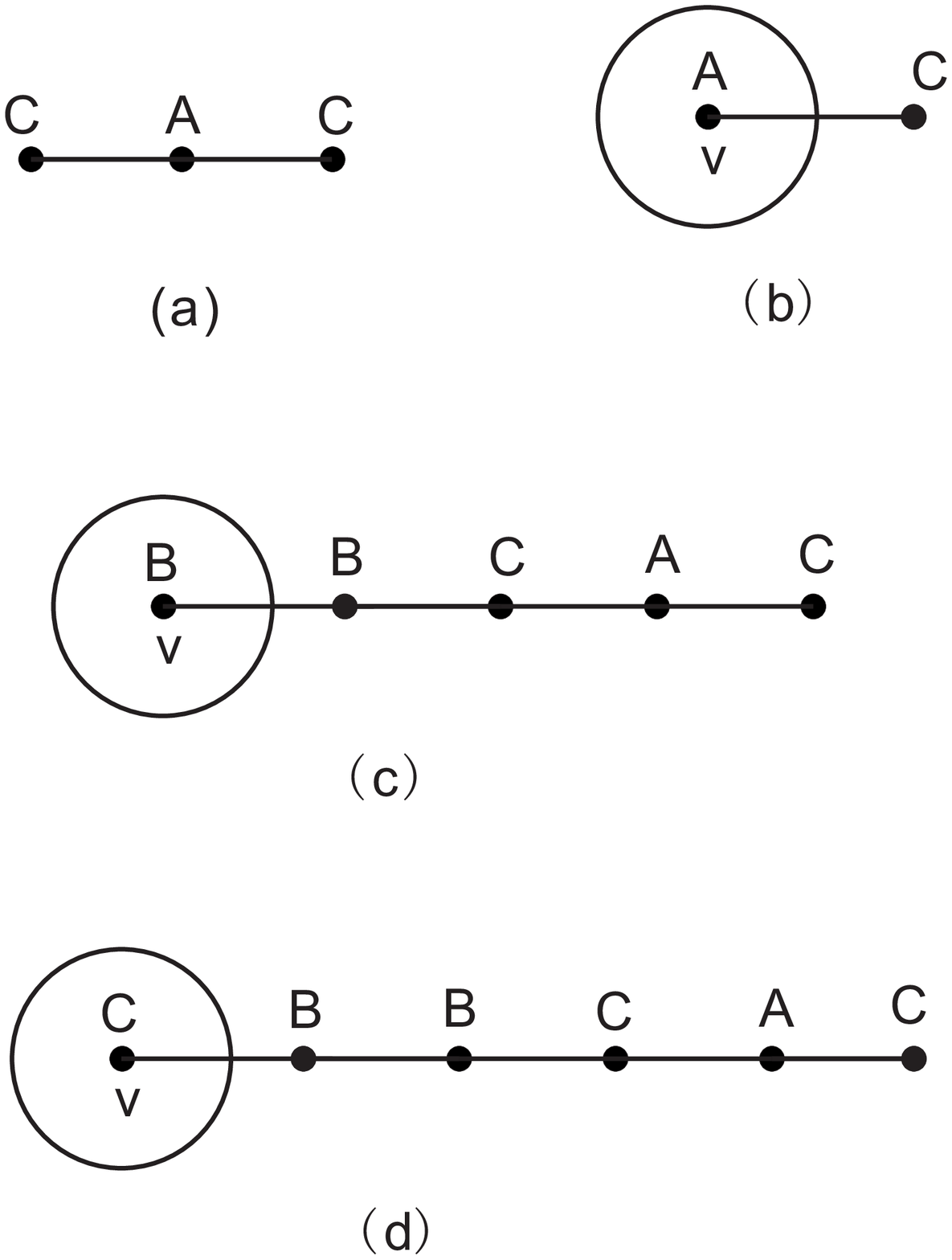}
 \end{center}
\qquad \qquad \qquad \qquad \qquad \qquad \qquad \qquad \qquad
{\small Fig.1} $\\$

Let $\mathscr{U}$ be the family of labeled trees that: (i) contains
$(P_3, S_0')$ where $S_0'$ is the labeling that assigns to the two
leaves of the path $P_3$ status $C$ and to the support vertex status
$A$ (see Fig.1(a)); and (ii) is closed under the operations
$\mathscr{P}_1$, $\mathscr{P}_2$ and $\mathscr{P}_3$ that are listed
below, which extend the tree $T'$ to a tree $T$ by attaching a tree
to the vertex $v\in V(T')$.

{\bf Operation} $\mathscr{P}_1$: Let $v$ be a vertex with
sta$(v)=A$. Add a vertex $u$ and the edge $uv$. Let sta$(u)=C$.

{\bf Operation} $\mathscr{P}_2$: Let $v$ be a vertex with
sta$(v)=B$. Add a path $v_1v_2v_3v_4$ and the edge $vv_1$. Let
sta$(v_1)=B$, sta$(v_2)=$sta$(v_4)=C$, and sta$(v_3)=A$.

{\bf Operation} $\mathscr{P}_3$: Let $v$ be a vertex with
sta$(v)=C$. Add a path $v_1v_2v_3v_4v_5$ and the edge $vv_1$. Let
sta$(v_1)=$sta$(v_2)=B$, sta$(v_3)=$sta$(v_5)=C$ and sta$(v_4)=A$.

$\\$

\begin{center}
  \includegraphics[width=3.2in]{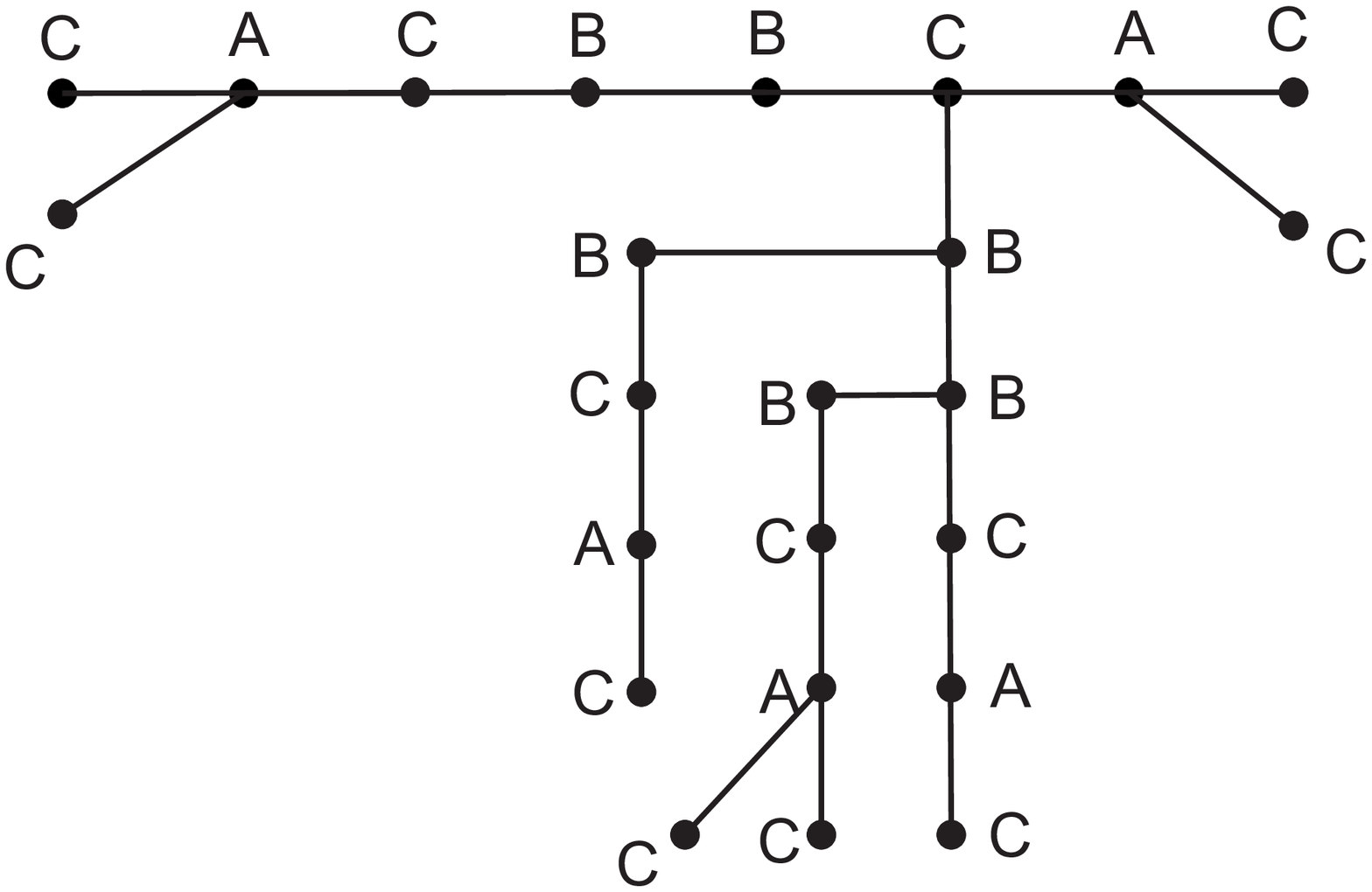}
 \end{center}
\qquad \qquad \qquad \qquad \qquad \qquad \qquad \qquad \qquad
{\small Fig.2} $\\$

Let $(T, S)\in \mathscr{U}$ be a labeled tree for some labeling $S$.
Then there is a sequence of labeled trees $(P_3, S_0')$, $(T_1, S_1),
\cdots, (T_{k-1}, S_{k-1})$, $(T_k, S_k)$ such that $(T_k, S_k)=(T,
S)$. The labeled tree $(T_i, S_i)$ can be obtained from $(T_{i-1},
S_{i-1})$ by one of the operations $\mathscr{P}_1$, $\mathscr{P}_2$ and
$\mathscr{P}_3$, where $i\in \{1, 2, \cdots, k\}$. We
remark that a sequence of labeled trees used to construct $(T, S)$
is not necessarily unique. The graph in Fig.2 is an example which belongs to $\mathscr{U}$.

In what follows, we present a few preliminary results.

\begin{obs}
Let $T$ be a tree of order at least $3$ and $S$ be a labeling of $T$
such that $(T, S)\in \mathscr{U}$. Then, $T$ has the following
properties:

$(a)$ If a vertex $v$ is a support vertex, then sta$(v)=A$, and each
of its neighbors has status $C$.

$(b)$ If a vertex $v$ is a leaf, then $v$ has status $C$.

$(c)$ If a vertex $v$ is labeled $C$, then it is a leaf or a vertex
all of whose neighbors are labeled $B$ except for one, which is
labeled $A$.

$(d)$ If a vertex $v$ has status $B$, then all of the neighbors of
$v$ are labeled $B$ except for one, which is labeled $C$.

$(e)$ $S_A$ and $S_C$ are two independent sets of $T$.
\end{obs}

Before giving the following lemma, we shall need an additional
notation. We call $D$ an \emph{almost semitotal dominating set} of a
graph $G$ relative to a vertex $v$ if $D$ is a dominating set of $G$
and every vertex in $D$ is within distance $2$ of another vertex of
$D$, except for $v$.

\begin{lem}
If $T$ is a tree such that $(T, S)\in \mathscr{U}$ for some labeling
$S$, then for any vertex $x\in S_A$, there exists an almost
semitotal dominating set of $T$ relative to $x$ with cardinality
$\gamma_{t2}(T)-1$.
\end{lem}

\bp We know that $(T, S)\in \mathscr{U}$ for some labeling $S$, and
as mentioned above, a sequence of labeled trees used to
construct $(T, S)$ is not necessarily unique. So we select a
sequence used to construct $(T, S)$: $(T_0, S_0)$, $(T_1, S_1),
\cdots, (T_{k-1}, S_{k-1})$, $(T_k, S_k)$, where $(T_0, S_0)=(P_3,
S_0')$ and $(T_k, S_k)=(T, S)$, such that the vertex $x\in
V(T_i)\setminus V(T_{i-1})$ ($i\in \{0, 1, 2, \cdots, k\}$) and the
number $i$ as small as possible (this condition is essential for the
following algorithm).

Now, we construct a set $H$ as follows.

\vspace{0.2cm}

$(I)$ Set $P:=\emptyset$ and $H:=\{t\}$, where $t$ is the vertex of
$V(T_0)$ which has status $A$ in $S_0$. Set $j:=1$.

\vspace{0.2cm}

$(II)$ We query whether $j>i$ or not.

--- If the answer to the query is `yes',

then go to $(IV)$.

--- If the answer to the query is `no',

then go to $(III)$.

\vspace{0.2cm}

$(III)$ We query which operation is used at the $j$-th step.

--- If $(T_j, S_j)$ is obtained from $(T_{j-1},
S_{j-1})$ by operation $\mathscr{P}_1$.

then set $j:=j+1$. Go to $(II)$.

--- If $(T_j, S_j)$ is obtained from $(T_{j-1},
S_{j-1})$ by operation $\mathscr{P}_2$.

then set $H:=H\cup P\cup \{y\}$, where $y$ is the vertex of
$V(T_j)\setminus V(T_{j-1})$ which has status $A$ in $S_j$. Set
$P:=\emptyset$ and put the vertex of $V(T_j)\setminus V(T_{j-1})$
which at distance $2$ from $y$ into $P$. Set $j:=j+1$. Go to $(II)$.

--- If $(T_j, S_j)$ is obtained from $(T_{j-1},
S_{j-1})$ by operation $\mathscr{P}_3$.

then set $H:=H\cup \{y, z\}$, where $y$ is the vertex of
$V(T_j)\setminus V(T_{j-1})$ which has status $A$ in $S_j$ and $z$
is the vertex of $V(T_j)\setminus V(T_{j-1})$ which at distance $3$
from $y$. Set $P:=\emptyset$ and put the vertex of $V(T_j)\setminus
V(T_{j-1})$ which at distance $2$ from $y$ into $P$. Set $j:=j+1$.
Go to $(II)$.

\vspace{0.2cm}

$(IV)$ We query whether $j>k$ or not.

--- If the answer to the query is `yes',

then we terminate.

--- If the answer to the query is `no',

then go to $(V)$.

 \vspace{0.2cm}

$(V)$ We query whether $|V(T_j)\setminus V(T_{j-1})|=1$ or not.

--- If the answer to the query is `yes',

then set $j:=j+1$. Go to $(IV)$.

--- If the answer to the query is `no',

then set $H:=H\cup \{w, h\}$, where $w$ is the vertex of
$V(T_j)\setminus V(T_{j-1})$ which has status $A$ in $S_j$, and $h$
is the vertex of $V(T_j)\setminus V(T_{j-1})$ which at distance $2$
from $w$. Set $j:=j+1$. Go to $(IV)$.

 \vspace{0.2cm}

After the end of this procedure, the set $H$ is a desire set.
Moreover, it follows from the method of constructing the set $H$
that $H$ contains all vertices of $S_A$. \ep

\begin{lem}
If $T$ is a tree such that $(T, S)\in \mathscr{U}$ for some labeling
$S$, then $T$ is in Class~3.
\end{lem}

\bp Let $T^{*}$ be obtained from $T$ by subdividing any edge $w$ of
$T$ twice. It is easy to see that $\gamma_{t2}(T)\leq
\gamma_{t2}(T^{*})$. In order to show that $T$ is in Class~3, we
need to show that $\gamma_{t2}(T)\geq \gamma_{t2}(T^{*})$.

Since $(T, S)\in \mathscr{U}$ for some labeling $S$, we can select a
sequence of labeled trees used to construct $(T, S)$: $(T_0, S_0)$,
$(T_1, S_1), \cdots, (T_{k-1}, S_{k-1})$, $(T_k, S_k)$, where $(T_0,
S_0)=(P_3, S_0')$ and $(T_k, S_k)=(T, S)$ such that $w\in
E(T_i)\setminus E(T_{i-1})$ and the number $i$ as small as possible.

If $i=0$, $w$ is a pendant edge in $T_0$, say $xx_1$, where $x$ is
the support vertex in $T_0$. Let $y_1, y_2$ be the two new vertices
resulting from subdividing the edge $xx_1$. By Lemma~2.8, there
exists an almost semitotal dominating set of $T$ relative to $x$
with cardinality $\gamma_{t2}(T)-1$, say $X$. Then, the set $X\cup
\{y_2\}$ is a semitotal dominating set of $T^{*}$. Hence,
$\gamma_{t2}(T^{*})\leq \gamma_{t2}(T)$. So we consider the case of
$i\neq 0$.

If $T_i$ is obtained from $T_{i-1}$ by adding a vertex $x_1$ and
joining it to a vertex $x_2$ of $T_{i-1}$, which has status $A$ in
$S_{i-1}$, then $w=x_1x_2$. We construct an almost semitotal
dominating set $H$ of $T$ relative to $x_2$ with cardinality
$\gamma_{t2}(T)-1$, the method of constructing the set $H$ is as
mentioned in the algorithm of Lemma~2.8. Let $y_1, y_2$ be the two
new vertices resulting from subdividing the edge $x_1x_2$. Then,
$H\cup \{y_1\}$ is a semitotal dominating set of $T^{*}$. That is,
$\gamma_{t2}(T^{*})\leq \gamma_{t2}(T)$.

If $T_i$ is obtained from $T_{i-1}$ by adding a path
$x_1x_2x_3x_4x_5$ and an edge $x_1x$, where $x$ has status $C$ in
$S_{i-1}$. We construct an almost semitotal dominating set $H$ of
$T$ relative to $x_4$ with cardinality $\gamma_{t2}(T)-1$, the
method of constructing the set $H$ is as mentioned in the algorithm
of Lemma~2.8. It follows from the construction method of $H$ and the
definition of almost semitotal dominating set that $x_1\in H$. Let
$y_1, y_2$ be the two new vertices resulting from subdividing the
edge $w$. If $w=xx_1$, then $(H\setminus \{x_1\})\cup \{y_1, x_2\}$
is a semitotal dominating set of $T^{*}$. If $w=x_1x_2$, then $H\cup
\{x_2\}$ is a semitotal dominating set of $T^{*}$. If $w\in
\{x_2x_3, x_3x_4, x_4x_5\}$, the proof is similar to the argument as
above. In either case, we have that $\gamma_{t2}(T^{*})\leq
\gamma_{t2}(T)$.

If $T_i$ is obtained from $T_{i-1}$ by adding a path $x_1x_2x_3x_4$
and an edge $x_1x$, where $x$ has status $B$ in $S_{i-1}$, the proof
is similar to the argument as above.
 \ep

 \begin{obs}
Let $G$ be a connected graph that is not a star. Then,

$(i)$ there is a $\gamma$-set that contains no leaf of $G$, and

$(ii)$\cite{Henning1} there is a $\gamma_{t2}$-set that contains no
leaf of $G$.
\end{obs}

\begin{lem}
If a tree $T$ of order at least $3$ is in Class~3, then $(T, S)\in
\mathscr{U}$ for some labeling $S$.
\end{lem}

\bp We proceed by induction on the order $n$ of $T$. If $T$ is a
star of order at least $3$, then it is in Class~3, and $(T, S)\in
\mathscr{U}$, where $S$ is the labeling that assigns to the support
vertex of $T$ status $A$ and to the leaves status $C$. It is easy to
verify that no tree whose diameter is at most $6$ is in Class~3,
except for the stars of order at least $3$. So we consider the case
that $diam(T)\geq 7$. Assume that for any tree $T'$ in Class~3 with
order less than $|T|$, we always have that $(T', S')\in \mathscr{U}$
for some labeling $S'$.

{\flushleft\textbf{Claim 1.}}\quad Each support vertex has exactly
one leaf-neighbor.

If not, assume that there is a support vertex $u$ which is adjacent
to at least two leaves. Deleting one of its leaf-neighbors, say
$u_1$, and denote the resulting tree by $T'$. Take an edge $w\in
E(T')$, let $T^{*}$ (respectively, $T'^{*}$) be obtained from $T$
(respectively, $T'$) by subdividing the edge $w$ twice. Let $D$ be a
$\gamma_{t2}$-set of $T'$ containing no leaf. Clearly, $D$ is a
semitotal dominating set of $T$. Then, we have that
$\gamma_{t2}(T)\leq \gamma_{t2}(T')\leq \gamma_{t2}(T'^{*})\leq
\gamma_{t2}(T^{*})=\gamma_{t2}(T)$. Thus we must have equality
throughout this inequality chain, whence $\gamma_{t2}(T')=
\gamma_{t2}(T'^{*})$. That is, $T'$ is in Class~3. By the inductive
hypothesis, $(T', S')\in \mathscr{U}$ for some labeling $S'$. Let
$S$ be obtained from the labeling $S'$ by labeling the vertex $u_1$
with label $C$. Then, $(T, S)$ can be obtained from $(T', S')$ by
operation $\mathscr{P}_1$. Thus, $(T, S)\in \mathscr{U}$.\ep

Let $P=v_1v_2v_3\cdots v_t$ be a longest path in $T$ such that

(i) $d(v_5)$ as large as possible, and subject to this condition

(ii) $d(v_4)$ as large as possible.

By Claim~1, $d(v_2)=2$. It follows from Observation~2.6 that
$d(v_3)=2$.

{\flushleft\textbf{Claim 2.}}\quad $d(v_4)=2$.

Assume that $d(v_4)>2$. From Observation~2.6(2), $v_4$ is not a
support vertex. Let $u$ be a neighbor of $v_4$ outside $P$. Then,
either $u$ is a support vertex of degree two, or $d(u)=2$ and it is
adjacent to a support vertex of degree two outside $P$.

In either case, we subdivide the edge $uv_4$ twice, and denote the
resulting tree by $T^{*}$. Clearly, $\gamma_{t2}(T^{*})-1\geq
\gamma_{t2}(T)$. Contradicting to the condition that $T$ is in
Class~3. \ep

{\flushleft\textbf{Claim 3.}}\quad $d(v_5)=2$.

Assume that $d(v_5)>2$. Let $u$ be a neighbor of $v_5$ outside $P$.
If $u$ is a leaf or a support vertex, we subdivide the edge $uv_5$
twice, and denote the resulting tree by $T^{*}$. Clearly,
$\gamma_{t2}(T^{*})-1\geq \gamma_{t2}(T)$. Contradicting to the
condition that $T$ is in Class~3.

Since $d(v_5)>2$, there exists the leaves outside $P$, say $a_1,
a_2, \cdots, a_l$, such that for each $i\in \{1, 2, \cdots, l\}$,
$V(P_i)\cap V(P)=\{v_5\}$, where $P_i$ is the shortest path between
$a_i$ and $v_5$. Without loss of generality, assume that
$P_1=v_5u_su_{s-1}\cdots u_1$ be the longest path among all $P_i$,
where $u_1=a_1$. Note that $s=3$ or $4$.

From Observation~2.6, Claim~1 and the choice of $P$, we only need to
consider the case that each $u_i$ has degree two, where $i=2, 3,
\cdots, s$.

Let $T'=T-\{v_1, v_2, v_3, v_4\}$. Clearly, $\gamma_{t2}(T)\leq
\gamma_{t2}(T')+2$. Let $T^{*}$ (respectively, $T'^{*}$) be obtained
from $T$ (respectively, $T'$) by subdividing an edge $w\in E(T')$
twice. Next, we ready to show that $\gamma_{t2}(T^{*})-2\geq
\gamma_{t2}(T'^{*})$. If $w\not \in \{v_5u_s, u_su_{s-1}, \cdots,
u_2u_1\}$, then we are done. So $w\in \{v_5u_s, u_su_{s-1}, \cdots,
u_2u_1\}$, without loss of generality, assume that $w=v_5u_s$. That
is, $T^{*}$ (respectively, $T'^{*}$) be obtained from $T$
(respectively, $T'$) by subdividing the edge $v_5u_s$ with vertices
$x_1, x_2$.

If $s=3$, let $H$ be obtained from $T$ by subdividing the edge
$v_5v_6$ with vertices $y_1, y_2$. Let $D$ be a $\gamma_{t2}(H)$-set
which contains no leaf. Then, $v_2\in D$. Note that $|\{v_3,
v_4\}\cap D|=1$. Without loss of generality, let $v_4\in D$ (If
$v_3$ belongs to $D$, then we can replace it in $D$ by $v_4$).
Similarly, we have that $u_2, v_5\in D$. Clearly, $|\{y_1, y_2\}\cap
D|\leq 1$. If $|\{y_1, y_2\}\cap D|=1$ and $v_6\not \in D$, we can
simply replace $x$ in $D$ by $v_6$, where $x\in \{y_1, y_2\}\cap D$.
If $|\{y_1, y_2\}\cap D|=1$ and $v_6\in D$, we can simply replace
$x$ in $D$ by $y$, where $x\in \{y_1, y_2\}\cap D$ and $y\in
N_H(v_6)\setminus \{y_2\}$ (Note that if $\{y_1, y_2\}\cap
D=\emptyset$, then $v_6\in D$). Let $D'=(D\setminus \{v_5\})\cup
\{x_2\}$ and $D''=(D\setminus \{v_2, v_4, v_5\})\cup \{x_2\}$. The
set $D'$ is a $\gamma_{t2}$-set of $T^{*}$, and the set $D''$ is a
semitotal dominating set of $T'^{*}$. That is,
$\gamma_{t2}(T^{*})-2\geq \gamma_{t2}(T'^{*})$.

If $s=4$, by a similar argument as above, we can also obtain the
same conclusion. That is, $\gamma_{t2}(T^{*})-2\geq
\gamma_{t2}(T'^{*})$.

In summary, we have that $\gamma_{t2}(T)\leq \gamma_{t2}(T')+2\leq
\gamma_{t2}(T'^{*})+2\leq \gamma_{t2}(T^{*})=\gamma_{t2}(T)$.
Consequently we must have equality throughout this inequality chain,
whence $\gamma_{t2}(T')=\gamma_{t2}(T'^{*})$. It follows that $T'$
is in Class~3. By induction, $(T', S')\in \mathscr{U}$ for some
labeling $S'$. And then, $u_1, u_2$ have status $C$ and $A$,
respectively. Moreover, by Observation~2.7(a), (c) and (d), $u_3,
v_5$ have status $C, B$ respectively when $s=3$, and $u_3, u_4, v_5$
have status $C, B, B$ respectively when $s=4$. In either case, let
$S$ be obtained from the labeling $S'$ by labeling the vertex $v_1,
v_2, v_3, v_4$ with label $C, A, C, B$, respectively. Then, $(T, S)$
can be obtained from $(T', S')$ by operation $\mathscr{P}_2$. Thus,
$(T, S)\in \mathscr{T}$.\ep

Now we let $T'=T-\{v_1, v_2, v_3, v_4, v_5\}$, let $w\in E(T')$ and
$T^{*}$ (respectively, $T'^{*}$) be obtained from $T$ (respectively,
$T'$) by subdividing the edge $w$ twice. Clearly, we have that
$\gamma_{t2}(T')+2\geq \gamma_{t2}(T)$.

If $d(v_6)>2$, then $v_6$ is not a support vertex. Otherwise, let
$H$ be obtained from $T$ by subdividing the edge $v_1v_2$ twice. It
is easy to verify that $\gamma_{t2}(H)-1\geq \gamma_{t2}(T)$.
Contradicting to the assumption that $T$ is in Class~3. Since
$d(v_6)>2$, there exists the leaves, say $b_1, b_2, \cdots, b_l$,
such that for each $i\in \{1, 2, \cdots, l\}$, $V(P_i')\cap
V(P)=\{v_6\}$, where $P_i'$ is the shortest path between $b_i$ and
$v_6$. Without loss of generality, assume that $P_1'=v_6u_0u_1\cdots
u_s$ be the longest path among all $P_i'$, where $u_s=b_1$. Note
that $s\leq 4$. We only need to consider the case that $d(u_i)=2$
for $i=0, 1, \cdots, s-1$ (otherwise, the proof is similar to the
previous arguments).

If $s=2$ or $3$, we can obtain a similar contradiction as above. So
we only need to consider the case that $d(v_6)=2$, or $d(v_6)>2$ and
$s=1, 4$. In these cases, by similar arguments as in Claim~3, we
have that $\gamma_{t2}(T^{*})-2\geq \gamma_{t2}(T'^{*})$. Hence,
$\gamma_{t2}(T)\leq \gamma_{t2}(T')+2\leq \gamma_{t2}(T'^{*})+2\leq
\gamma_{t2}(T^{*})=\gamma_{t2}(T)$. Thus we must have equality
throughout this inequality chain, whence $\gamma_{t2}(T')=
\gamma_{t2}(T'^{*})$. That is, $T'$ is in Class~3. By the inductive
hypothesis, $(T', S')\in \mathscr{U}$ for some labeling $S'$.

If $d(v_6)=2$, or $d(v_6)>2$ and $s=1$, then $v_6$ has status $C$ in
$S'$. If $d(v_6)>2$ and $s=4$, the vertices $u_4, u_3, u_2, u_1,
u_0$ have status $C, A, C, B, B$ in $S'$, respectively. Since
$d(u_0)=2$, by Observation~2.7(d), $v_6$ has status $C$. In either
case, let $S$ be obtained from the labeling $S'$ by labeling the
vertices $v_1, v_2, v_3, v_4, v_5$ with label $C, A, C, B, B$,
respectively. Then, $(T, S)$ can be obtained from $(T', S')$ by
operation $\mathscr{P}_3$. Thus, $(T, S)\in \mathscr{U}$.\ep

As an immediate consequence of Lemmas~2.9 and 2.11 we have the
following conclusion.

\begin{thm}
A tree $T$ of order at least $3$ is in Class~3 if and only if $(T,
S)\in \mathscr{U}$ for some labeling $S$.
\end{thm}

\section{Semitotal domination versus domination and total domination in trees}

\subsection{Main results}

It is well known that the
semitotal domination number is a parameter that is squeezed between domination number and total domination number, namely $\gamma(G)\leq \gamma_{t2}(G)\leq \gamma_t(G)$. So it is natural to consider the ratios $\frac{\gamma_{t2}(T)}{\gamma(T)}$ and $\frac{\gamma_t(T)}{\gamma_{t2}(T)}$.

In what follows, we give three operations as follows:

{\bf Operation} $\mathscr{O}_1$: Let $v$ be a vertex with sta$(v)=A
$ or $B$. Add a vertex $u$ and the edge $uv$. Let sta$(u)=C$.

{\bf Operation} $\mathscr{O}_2$: Let $v$ be a vertex with
sta$(v)=A$. Add a path $u_1u_2u_3u_4$ and the edge $u_1v$. Let
sta$(u_1)=D$, sta$(u_2)=E$, sta$(u_3)=B$ and sta$(u_4)=C$.

{\bf Operation} $\mathscr{O}_3$: Let $v$ be a vertex with
sta$(v)=A$. Add a path $u_1u_2u_3$ and the edge $u_1v$. Let
sta$(u_1)=D$, sta$(u_2)=B$, sta$(u_3)=C$.

The three operations $\mathscr{O}_1$, $\mathscr{O}_2$ and $\mathscr{O}_3$ are
illustrated in Fig.3(c), (d) and (e).

Next, we are ready to give two families $\mathscr{T}$ and $\mathscr{T}_1$.

Let $\mathscr{T}$ be the family of labeled trees that: (i) contains
$(P_6, S_0)$ where $S_0$ is the labeling that assigns to the two
leaves of the path $P_6$ status $C$, to the two support vertices
status $A$ and $B$ respectively, and to the two center vertices
status $D$ and $E$ respectively (see Fig.3(a)); and (ii) is closed
under the two operations $\mathscr{O}_1$ and $\mathscr{O}_2$ that
are listed below, which extend the tree $T'$ to a tree $T$ by
attaching a tree to the vertex $v\in V(T')$.

Let $\mathscr{T}_1$ be the family of labeled trees that: (i) contains
$(P_5, S_0')$ where $S_0'$ is the labeling that assigns to the two
leaves of the path $P_5$ status $C$, to the two support vertices
status $B$ and $A$ respectively, and to the center vertex status
$D$  (see Fig.3(b)); and (ii) is closed under the two operations $\mathscr{O}_1$ and
$\mathscr{O}_3$ that are listed below, which extend the tree $T'$ to
a tree $T$ by attaching a tree to the vertex $v\in V(T')$.

$\\$

\begin{center}
  \includegraphics[width=4.5in]{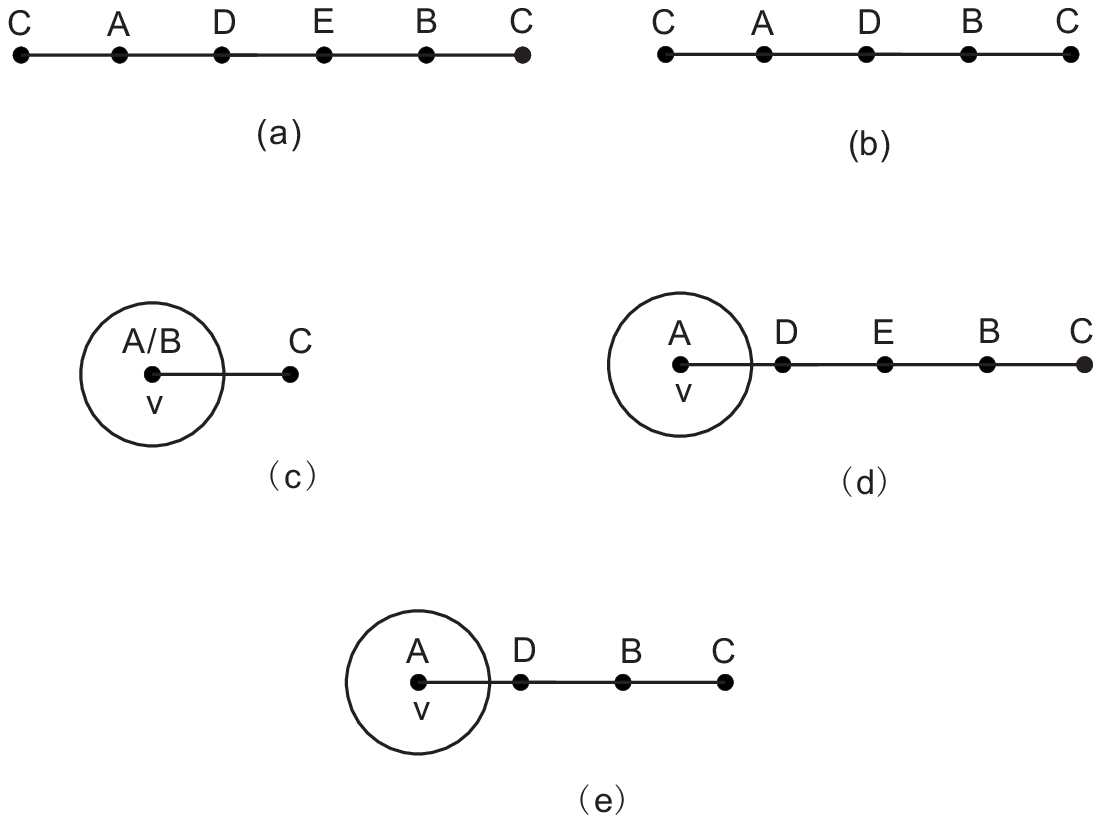}
 \end{center}
\qquad \qquad \qquad \qquad \qquad \qquad \qquad \qquad \qquad
{\small Fig.3} $\\$

Two main conclusions of this section are listed as follows.

\begin{thm}
Let $T$ be a tree that is not a star, we have that
$\gamma_{t2}(T)\leq 2\gamma(T)-1$. Moreover, the trees $T$
satisfying $\gamma_{t2}(T)=2\gamma(T)-1$ are precisely those trees
$T$ such that $(T, S)\in \mathscr{T}$ for some labeling $S$.
\end{thm}

\begin{thm}
Let $T$ be a nontrivial tree, we have that $\gamma_t(T)\leq
2\gamma_{t2}(T)-1$. Moreover, the trees $T$ satisfying
$\gamma_t(T)=2\gamma_{t2}(T)-1$ are precisely those trees $T$ such
that $(T, S')\in \mathscr{T}_1$ for some labeling $S'$.
\end{thm}

The proofs of the above two theorems are very similar. Due to the readability of this article, we will only give the proof of Theorem~3.1.

\subsection{The proof of Theorem~3.1}

In what follows, we present a few preliminary results.

\begin{obs}
Let $T$ be a tree of order at least $6$ and $S$ be a labeling of $T$
such that $(T, S)\in \mathscr{T}$. Then, $T$ has the following
properties:

$(a)$ A vertex is labeled $A$ or $B$ if and only if it is a support
vertex.

$(b)$ A vertex is labeled $C$ if and only if it is a leaf.

$(c)$ $|S_A|=1$, $|S_B|=|S_D|=|S_E|$.

$(d)$ The set $S_A\cup S_B$ is the unique $\gamma$-set of $T$.

$(e)$ The set $S_A\cup S_B\cup S_D$ is a $\gamma_{t2}$-set of $T$.

$(f)$ If a vertex has status $A$ (respectively, $B$), then each of
its non-leaf neighbors is labeled $D$ (respectively, $E$).

$(g)$ If a vertex has status $D$ (respectively, $E$), then it has
degree two and the two neighbors are labeled $A$ and $E$
(respectively, $B$ and $D$).
\end{obs}

From Observation~3.3 (c), (d) and (e), the following corollary can
be derived immediately.

\begin{cor}
Let $T$ be a tree and $S$ be a labeling of $T$ such that $(T, S)\in
\mathscr{T}$. Then, $\gamma_{t2}(T)=2\gamma(T)-1$.
\end{cor}

Next, we give the proof of Theorem~3.1.

\bp The sufficiency follows immediately from Corollary~3.4. So we prove the necessity only.
The proof is by induction on the order of $T$. If $|T|\leq 6$, it
is easy to verify that $\gamma_{t2}(T)\leq 2\gamma(T)-1$, and
$T=P_6$ when the equality holds. So we let $|T|\geq 7$ and assume
that for every non-star tree $T'$ of order less than $|T|$ we have
$\gamma_{t2}(T')\leq 2\gamma(T')-1$, with equality if and only if
$(T', S')\in \mathscr{T}$ for some labeling $S'$. By
Observation~2.10(i), there exists a $\gamma$-set of $T$ which
contains no leaf, say $D$.

{\flushleft\textbf{Claim 1.}}\quad Each support vertex has exactly
one leaf-neighbor.

Suppose that $v$ is a support vertex which has at least two
leaf-neighbors, say $v_1, v_2$. Let $T'=T-v_1$ and $R$ be a
$\gamma_{t2}$-set of $T'$ containing no leaf. Then, $R$ is also a
semitotal dominating set of $T$. Hence, $\gamma_{t2}(T)\leq
\gamma_{t2}(T')$. Combining the fact that $\gamma(T')\leq
\gamma(T)$, we have that $\gamma_{t2}(T)\leq \gamma_{t2}(T')\leq
2\gamma(T')-1\leq 2\gamma(T)-1$. If $\gamma_{t2}(T)=2\gamma(T)-1$,
then we have that $\gamma_{t2}(T')=2\gamma(T')-1$. By the inductive
hypothesis, $(T', S')\in \mathscr{T}$ for some labeling $S'$. It
follows from Observation~3.3(a) that $v$ has status $A$ or $B$ in
$S'$. Let $S$ be obtained from $S'$ by labeling the vertex $v_1$
with label $C$. Then, $(T, S)$ can be obtained from $(T', S')$ by
operation $\mathscr{O}_1$. Thus, $(T, S)\in \mathscr{T}$.\ep

We suppose that $diam(T)\geq 6$ (the result is trivial when
$diam(T)\leq 5$) and $P=v_1v_2v_3\cdots v_t$ be a longest path in
$T$ such that $d(v_3)$ as large as possible.

{\flushleft\textbf{Claim 2.}}\quad $d(v_3)=2$.

Assume that $d(v_3)>2$. Let $T'=T-\{v_1, v_2\}$ and $R'$ be a
$\gamma_{t2}$-set of $T'$ containing no leaf. Note that $D\setminus
\{v_2\}$ is a dominating set of $T'$. On the other hand, $R'\cup
\{v_2\}$ be a semitotal dominating set of $T$. Thus,
$\gamma_{t2}(T)\leq \gamma_{t2}(T')+1\leq 2\gamma(T')-1+1\leq
2\gamma(T)-2$.\ep

{\flushleft\textbf{Claim 3.}}\quad $d(v_4)=2$.

Assume that $d(v_4)>2$ and $u_1$ is a neighbor of $v_4$ outside $P$.
Let $T'=T-\{v_1, v_2, v_3\}$ and $R'$ be a $\gamma_{t2}$-set of $T'$
containing no leaf. From Claim~1 and the choice of $P$, we have that
at least one of the two conditions as follows holds:

(1) $u_1$ is a leaf;

(2) $u_1$ is a support vertex of degree two;

(3) $u_1$ has degree two and is adjacent to a support vertex of
degree two, say $u_2$.

In the first case, $v_4$ belongs to $D$ and $R'$. Hence,
$\gamma(T)-1\geq \gamma(T')$ and $\gamma_{t2}(T')+1\geq
\gamma_{t2}(T)$. Similar to the proof of Claim~2, we have that
$\gamma_{t2}(T)\leq 2\gamma(T)-2$.

In the second case, $u_1$ belongs to $D$ and $R'$. Then,
$\gamma(T)-1\geq \gamma(T')$ and $\gamma_{t2}(T')+2\geq
\gamma_{t2}(T)$. It means that $\gamma_{t2}(T)\leq
\gamma_{t2}(T')+2\leq 2\gamma(T')-1+2\leq 2\gamma(T)-1$. Suppose
next that $\gamma_{t2}(T)=2\gamma(T)-1$. Then we have equality
throughout the above inequality chain. In particular,
$\gamma_{t2}(T')=2\gamma(T')-1$. By induction, $(T', S')\in
\mathscr{T}$ for some labeling $S'$. Since $u_1$ is a support vertex
in $T'$, it follows from Observation~3.3(a) that $u_1$ has status
$A$ or $B$ in $S'$.

If sta$(u_1)=A$, from $d(u_1)=2$ and Claim~1, we have that $T'=P_6$
and $v_4$ has status $D$ in $S'$. Then, $T$ is the tree obtained
from a star of order four by subdividing two edges twice and the
remaining edge once. But in this case, $\gamma_{t2}(T)=4$ and
$\gamma(T)=3$, a contradiction.

If sta$(u_1)=B$, then $v_4$ has status $E$ in $S'$. It is easy to
check that $\gamma_{t2}(T)=2\gamma(T)-2$, a contradiction.

In the third case, $u_2$ belongs to $D$ and $R'$, $|\{u_1, v_4\}\cap
R'|=1$. Without loss of generality, we let $v_4\in R'$ (If $u_1\in
R'$, then we can replace $u_1$ in $R'$ by $v_4$). It follows that
$\gamma(T)-1\geq \gamma(T')$ and $\gamma_{t2}(T')+1\geq
\gamma_{t2}(T)$. It means that $\gamma_{t2}(T)\leq 2\gamma(T)-2$.\ep

Now, we let $T'=T-\{v_1, v_2, v_3, v_4\}$ and $R'$ be a
$\gamma_{t2}$-set of $T'$. Clearly, $|\{v_3, v_4, v_5\}\cap D|=1$.
Without loss of generality, $v_5\in D$ (If $v_3$ or $v_4$ belongs to
$D$, then we can replace it in $D$ by $v_5$). It implies that
$\gamma(T)-1\geq \gamma(T')$. On the other hand, $R'\cup \{v_2,
v_3\}$ be a semitotal dominating set of $T$. Hence,
$\gamma_{t2}(T)\leq \gamma_{t2}(T')+2\leq 2\gamma(T')-1+2\leq
2\gamma(T)-1$. Suppose next that $\gamma_{t2}(T)=2\gamma(T)-1$. Then
we have equality throughout the above inequality chain. In
particular, $\gamma_{t2}(T')=2\gamma(T')-1$. By induction, $(T',
S')\in \mathscr{T}$ for some labeling $S'$.

If $v_5$ has status $A$ in $S'$, let $S$ be obtained from $S'$ by
labeling the vertices $v_1, v_2, v_3, v_4$ with label $C, B, E, D$,
respectively. Then, $(T, S)$ can be obtained from $(T', S')$ by
operation $\mathscr{O}_2$. Thus, $(T, S)\in \mathscr{T}$. If
sta$(v_5)\in \{C, D, E\}$ in $S'$, or sta$(v_5)=B$ in $S'$ and
$|T'|>6$, we always have that $\gamma_{t2}(T)\leq 2\gamma(T)-2$, a
contradiction. If sta$(v_5)=B$ in $S'$ and $|T'|=6$, then $T'=P_6$.
Moreover, $v_6, v_7, v_8$ have status $E, D, A$ in $S'$,
respectively. Let $S''$ be obtained from $S'$ by relabeling the
vertices $v_5, v_6, v_7, v_8$ with label $A, D, E, B$, respectively.
And let $S$ be obtained from $S''$ by labeling the vertices $v_1,
v_2, v_3, v_4$ with label $C, B, E, D$, respectively. Thus, we can
also obtain that $(T, S)\in \mathscr{T}$.\ep


\end{document}